# ESTIMATION OF DISTRIBUTIONS, MOMENTS AND QUANTILES IN DECONVOLUTION PROBLEMS

BY PETER HALL AND SOUMENDRA N. LAHIRI

*University of Melbourne and University of Melbourne and Texas A&M University*

When using the bootstrap in the presence of measurement error, we must first estimate the target distribution function; we cannot directly resample, since we do not have a sample from the target. These and other considerations motivate the development of estimators of distributions, and of related quantities such as moments and quantiles, in errors-in-variables settings. We show that such estimators have curious and unexpected properties. For example, if the distributions of the variable of interest, $W$, say, and of the observation error are both centered at zero, then the rate of convergence of an estimator of the distribution function of $W$ can be slower at the origin than away from the origin. This is an intrinsic characteristic of the problem, not a quirk of particular estimators; the property holds true for optimal estimators.

**1. Introduction.** The problem of nonparametrically estimating a probability density, when the data are observed with error, has attracted a great deal of interest. However, in a range of circumstances the practical implementation of such estimators can be unattractive, since convergence rates are slow. Moreover, it is the distribution function, and not the density, that is needed in a wide variety of settings. For example, while in conventional applications of the bootstrap we proceed by resampling, and do not need to compute an empirical distribution function, this approach is infeasible when measurement errors are present; instead, we must generate data via a distribution-function estimate.

Therefore, in measurement-error problems, explicit distribution-function estimation assumes a substantial degree of importance which it does not necessarily enjoy in other settings. However, distribution-function estimators enjoy properties very different from those of their density counterparts. In particular, root-$n$ consistent estimation of a distribution function is possible if the error distribution is not too smooth. We shall give a necessary and sufficient condition for there to exist distribution-function estimators that converge at rate $n^{-1/2}$, and we shall explore, both theoretically and numerically, their intriguing properties. For example, we shall show that faster convergence rates can be achieved away from the origin than close to the origin.

It would be misleading to treat this problem in isolation; the unusual properties of distribution estimators are reflected in estimators of smooth functionals of distributions, for example, in quantile and moment estimators. However, while estimators in both these settings can be root-$n$ consistent, unusual features make the problems intrinsically interesting. In particular, while any polynomial moment can be estimated root-$n$ consistently, where $n$ denotes sample size, this is not true of fractional moments. In such cases, root-$n$ consistency is feasible if and only if the error distribution is not smoother than a certain amount, where the latter condition becomes less stringent as the exponent of the moment increases. When root-$n$ consistency is possible, it can be achieved without any statistical smoothing. In other cases, however, smoothing is necessary in order to achieve minimax-optimal convergence rates.

To give a little background to the problem of distribution estimation, we mention that, toward the end of his seminal paper on deconvolution density estimation, Fan (1991a) explored the distribution-estimation problem. He noted that upper bounds, for his particular estimator, and minimax lower bounds for arbitrary estimators, could be obtained, but found that they were of different orders of magnitude. He conjectured that his upper bound gave the optimal rate, and that the lower-bound rate could be increased to that of the upper bound. He suggested that the reason for the gap might be that the problem is more complex than his two-alternative analysis allowed, and that a highly composite-alternative approach could be necessary, as used by Stone (1982) in a different problem.

In fact, the problem is both simpler and more complex than this. It is simpler in the sense that a composite-alternative approach is not necessary in order to derive optimal rates, but more complex from the viewpoint that, apparently unsuspected by previous workers, the distribution-function estimator converges in an uneven fashion. Specifically, if the distributions of the variable of interest, $W$, say, and of the observation error are both centered at zero, then the rate of convergence of an estimator of the distribution function of $W$ can be relatively slow near the origin, with the result that the rate of convergence uniformly on the real line is an order of magnitude slower than the rate in the region $\{x : |x| > x_0\}$, for each fixed $x_0 > 0$. This remark applies both to the upper bound, for a particular estimator based on integrating a density estimator, and the lower bound, for arbitrary estimators. Therefore, uneven convergence rates are intrinsic to the problem, and are not artifacts of either our methodology or our mathematical arguments for deriving upper bounds.

Fan's (1991a) rates are in a slightly different context from ours; he measures smoothness in terms of derivatives, whereas we frame it through tail behavior of characteristic functions. The latter approach is arguably more natural in the present setting, because popular estimators are based on Fourier inversion. However, the two approaches can be reconciled closely. To the extent that this is possible, and in the context described in the previous paragraph, Fan's lower bound gives the optimal rate at the origin, although not at other places, while his upper bound is a little larger.

The context of density estimation has received greatest attention in the literature. Early contributions to this topic, suggesting estimators and discussing accuracy, include those of Carroll and Hall (1988), Devroye (1989), Stefanski and Carroll (1990), Zhang (1990) and Fan (1991a, 1991b, 1993). Hesse (1999) and Delaigle and Gijbels (2004a, 2004b) proposed methods for smoothing-parameter choice, Koo (1999) introduced a logspline-based deconvolution density estimator, Delaigle and Gijbels (2002) and Hesse and Meister (2004) discussed methods for estimating density derivatives, and van Es, Spreij and van Zanten (2003) treated volatility density estimation.

Recent contributions to optimality theory, in the context of density estimation, include those of Butucea (2004), who gave minimax convergence rates in cases where the unknown density belongs to a class of supersmooth functions, and the error distribution is ordinary-smooth; and Butucea and Tsybakov (2008), who provided sharp optimality results in settings where the unknown density and unknown error distribution are both supersmooth. Although practitioners have demonstrated a marked preference for kernel methods, wavelet-based deconvolution density estimators have been shown to enjoy excellent adaptivity properties. Note, for example, the contributions of Pensky and Vidakovic (1999), Fan and Koo (2002) and Pensky (2002), who derived convergence rates.

Groeneboom and Jongbloed (2003) discussed density estimators based on nonparametric maximum likelihood estimation when the error has a uniform distribution; see also Groeneboom and Wellner (1992). In terms of convergence rates, our work is more nearly related to these contributions than to most others in the setting of density estimation.

More closely related still are the papers of Booth and Hall (1993), who treated interval estimation in errors-in-variables models; Hesse (1995), who gave upper bounds to convergence rates of deconvolution distribution estimators; van de Geer (1995), who addressed estimation of a linear integral functional in a mixture model; Cordy and Thomas (1997), who discussed nonparametric estimation of a distribution function when it can be modeled as a mixture; Jongbloed (1998), who studied isotonic estimation of a distribution function; Ioannides and Papanastassiou (2001), who treated distribution estimation in the case of dependent data; and Qin and Feng (2003) and Cui (2005), who developed asymptotic properties of estimators of known functions of the mean of the target distribution.

## 2. Methodology.

2.1. *Estimators $\hat{f}_W$ and $\widehat{F}_W$.* Assume we observe $X_j = W_j + \delta_j$ for $1 \leq j \leq n$, where the $W_j$'s and $\delta_j$'s are independent. If the density $f_\delta$ of $\delta$ is not known, then the density $f_W$ of $W$ is not identifiable from the $X_j$'s alone. Therefore, it is very common (see, e.g., the literature cited in Section 1) to assume a form for $f_\delta$. Only in cases where, for instance, additional data are available directly on $\delta$ [Diggle and Hall (1993) and Neumann (1997)], or replicated data are available on $X$, would this assumption be unnecessary.

A conventional estimator of the density $f_W$ of $W$ is given by

$$\hat{f}_W(x) = \hat{f}_W(x \mid h) = \frac{1}{nh} \Re \sum_{j=1}^{n} L\left(\frac{x - X_j}{h}\right), \tag{2.1}$$

where $\Re$ denotes real part,

$$L(u) = \frac{1}{2\pi} \int_{-\infty}^{\infty} e^{-itu} \frac{K^{\mathrm{Ft}}(t)}{f_\delta^{\mathrm{Ft}}(t/h)} \, dt, \tag{2.2}$$

$K$ is a kernel function (in particular, a function that integrates to 1), $K^{\mathrm{Ft}}(t) = \int e^{itx} K(x) \, dx$ is its Fourier transform, and $h > 0$ is a smoothing parameter. Note that $\hat{f}_W$ is well defined even if $f_W$ does not exist. Here and below we use the notation $f_\delta^{\mathrm{Ft}}$ and $f_W^{\mathrm{Ft}}$, for the characteristic functions of the distributions $F_\delta$ and $F_W$, without necessarily requiring the existence of the respective densities $f_\delta$ or $f_W$.

Under the common assumption that $K^{\mathrm{Ft}}$ is compactly supported and $f_\delta^{\mathrm{Ft}}$ does not vanish on the real line, the integral at (2.2) is well defined and finite. There is no loss of generality in assuming $K$ is symmetric, and seldom any loss in supposing the same for $f_\delta$. We shall make these simplifying assumptions below; they are almost invariably satisfied in practice. Then, $L$ is real-valued, and so the symbol $\Re$ may be dropped from (2.1).

The estimator $\widehat{F}_W$ is defined as simply the integral of $\hat{f}_W$ over $(-\infty, x]$, even in cases where $f_W$ does not exist. Details concerning its computation and interpretation, especially in the case $h = 0$, will be given in Appendix A.1.

2.2. *Moment estimators.* If we wished to estimate a moment of $W$, say, $\mu_r = E(W^r)$, where $r \geq 1$ was an integer, a naive approach would be to base the estimator directly on empirical moments of $X$ and the known theoretical moments of $\delta$. Since symmetry of $F_\delta$ implies $E(\delta) = 0$, then

$$\mu_r = E(X^r) - \sum_{j=2}^{r} \binom{r}{j} E(\delta^j) \mu_{r-j}. \tag{2.3}$$

[Of course, $E(\delta^j)$ vanishes for odd $j$.] Given estimators $\tilde{\mu}_j$ of $\mu_j$ for $1 \leq j \leq r-2$, substitution into (2.3) suggests an estimator of $\mu_r$, for $r \geq 1$:

$$\tilde{\mu}_r = \frac{1}{n} \sum_{j=1}^{n} X_j^r - \sum_{j=2}^{r} \binom{r}{j} E(\delta^j) \tilde{\mu}_{r-j}. \tag{2.4}$$

In particular, $\tilde{\mu}_1 = \bar{X} = n^{-1} \sum_j X_i$ and $\tilde{\mu}_2 = n^{-1} \sum_j X_j^2 - E(\delta^2)$.

Exactly the same estimators are obtained using the empirical distribution function $\widehat{F}_W(\cdot \mid 0)$. That is, if we define

$$\widehat{\mu}_r = \lim_{h \to 0} \int_{-\infty}^{\infty} u^q \, d\widehat{F}_W(u \mid h), \tag{2.5}$$

then $\widehat{\mu}_r = \tilde{\mu}_r$ for $r \geq 1$.

Provided $E(W^{2r}) + E(\delta^{2r}) < \infty$, the estimator $\tilde{\mu}_r$, and hence also $\widehat{\mu}_r$, is root-$n$ consistent for $\mu_r$. However, root-$n$ consistency is generally not possible for estimators of absolute moments, such as $\nu_q = E|W|^q$, when $q > 0$ is not equal to a positive integer. There is no simple analogue of the estimator at (2.4) in this case, although $\widehat{\mu}_r$, at (2.5), is readily generalized to

$$\hat{\nu}_q(h) = \int_{-\infty}^{\infty} |u|^r \, d\widehat{F}_W(u \mid h).$$

We shall argue in Section 3.4 that, if $q$ is not an even integer, then $\hat{\nu}_q$ is root-$n$ consistent for $\nu_q$ if and only if $F_\delta$ is sufficiently "rough," expressed, for example, in terms of the rate of convergence of $f_\delta^{\text{Ft}}$ to zero in its tails. This condition is unnecessary when $q$ is an even integer.

2.3. *Quantile estimators.* To estimate the $u$th quantile, say, $\xi_u = F_W^{-1}(u)$, where $0 < u < 1$, we first render $\widehat{F}_W$ monotone by defining

$$\widehat{F}_W^{\text{mon}}(x) = \widehat{F}_W^{\text{mon}}(x \mid h) = \sup\{\widehat{F}_W(y \mid h) : y \leq x\},$$

and then we put

$$\hat{\xi}_u = \hat{\xi}_u(h) = (\widehat{F}_W^{\text{mon}})^{-1}(u) = \sup\{y : \widehat{F}_W^{\text{mon}}(y) \leq u\} = \sup\{y : \widehat{F}_W(y) \leq u\}.$$

Then, $\hat{\xi}_u$ is our estimator of $\xi_u$.

The monotonization step serves only to ensure that, with probability 1, $\hat{\xi}_u$ is well defined. For the choices of bandwidth, and values of $u$, that we use when establishing properties of $\hat{\xi}_u$, the probability that the monotonization step makes no difference to the value of $\hat{\xi}_u$, and, in particular, that $\hat{\xi}_u$ is well defined without it, converges to 1 as $n$ increases. In general, the mean-square convergence rate of $\hat{\xi}_u$ is strictly slower than $n^{-1}$, and depends on choice of $h$.

## 3. Theory related to optimality.

3.1. *Function classes.* Classes of functions indicated by $\mathcal{F}_j$ will be sets of distributions, $F_\delta$, say, of the error random variable $\delta$, while classes denoted by $\mathcal{G}_j$ will be sets of distributions, $F_W$, of $W$. The positive numbers $\alpha$ and $\beta$ will represent bounds to the degrees of the polynomial rates at which $f_\delta^{\mathrm{Ft}}(t)^{-1}$ and $f_W^{\mathrm{Ft}}(t)^{-1}$ diverge as $|t|$ increases. They are generally upper bounds in the case of $\alpha$, and lower bounds in that of $\beta$.

Given $C > 0$, write $\mathcal{F}_1(C)$ for the class of all distributions $F_\delta$ for which $f_\delta^{\mathrm{Ft}}$ is real-valued and positive everywhere, and

$$(3.1) \qquad \int_0^\infty t^{-2} \{f_\delta^{\mathrm{Ft}}(t)^{-1} - 1\}^2 \, dt \leq C.$$

The integral above is clearly finite on any compact set $[0, t_0]$, with $t_0 > 0$, and so (3.1) amounts to a condition on the rate at which the tails of $f_\delta^{\mathrm{Ft}}$ approach zero as $|t|$ increases. In particular, (3.1) can be viewed as holding if and only if $f_\delta^{\mathrm{Ft}}(t)$ does not converge to zero too quickly as $|t|$ increases, or equivalently, if and only if $F_\delta$ is not too smooth.

For example, $F_\delta \in \mathcal{F}_1(C)$ for sufficiently large $C > 0$, if the characteristic function satisfies

$$(3.2) \qquad f_\delta^{\mathrm{Ft}}(t) \geq B(1 + |t|)^{-\alpha}$$

for some $0 < \alpha < \frac{1}{2}$ and sufficiently large $B > 0$. Condition (3.2) is close to asserting that $F_\delta$ has at most $\alpha$ bounded derivatives. A symmetrized Gamma distribution, with density

$$(3.3) \qquad \begin{aligned} \phi_\alpha(x) &= \int_{|y|<\infty} \psi_\alpha(x+y)\psi_\alpha(y)\,dy \\ &\qquad \text{for } -\infty < x < \infty, \text{ where } \psi_\alpha(x) \\ &= \Gamma(\alpha/2)^{-1} x^{(\alpha/2)-1} e^{-x} \\ &\qquad \text{for } 0 < x < \infty \text{ and } \alpha > 0, \end{aligned}$$

satisfies both (3.1) and (3.2) provided $\alpha < \frac{1}{2}$.

Write $\mathcal{F}_2(C)$ for the class of all $F_\delta \in \mathcal{F}_1(C)$ for which $E|\delta| \leq C$. The function classes $\mathcal{F}_3(\alpha, C)$, $\mathcal{F}_4(C, q)$, $\mathcal{F}_5(\alpha, C)$ and $\mathcal{F}_6(\alpha, C)$ will be defined concisely in Appendix A.2. In heuristic terms, $\mathcal{F}_3(\alpha, C)$ is the class of distributions $F_\delta$ that satisfy (3.2), have a bounded density and bounded first absolute moment; and $\mathcal{F}_4(C, q)$ is a class of $F_\delta$ having sufficiently many finite moments and for which (3.1) holds but with the integral taken over $[1, \infty)$ and $t^{-2}$ replaced by $t^{-2(q+1)}$. The latter constraint increases the permitted smoothness of $F_\delta$, since it allows the tails of $f_\delta^{\mathrm{Ft}}$ to decrease relatively quickly.

The function class $\mathcal{F}_5(\alpha, C)$ is the set of distributions in $\mathcal{F}_3(\alpha, C)$ that have sufficiently many finite moments. And $\mathcal{F}_6(\alpha, C)$ is the subset of $\mathcal{F}_5(\alpha, C)$ for which the smoothness conditions on $f_\delta^{\mathrm{Ft}}$ are imposed not just on that function but, in an analogous way, on its first two derivatives as well.

For $C > 0$, let $\mathcal{G}_1(C)$ be the class of distributions $F_W$ that have densities $f_W$ satisfying $\sup_w f_W(w) \leq C$. Write $\mathcal{G}_2(C)$ for the class of $F_W$ for which $E|W| \leq C$. Note that there are distributions in the class $\mathcal{G}_2(C)$ for which $f_W$ does not exist.

The function classes $\mathcal{G}_3(\beta, C)$, $\mathcal{G}_4(\beta, C)$, $\mathcal{G}_5(C, k)$, $\mathcal{G}_6(\beta, C)$, $\mathcal{G}_7(\beta, C, u, g)$ and $\mathcal{G}_8(\beta, C, u, g)$ will be detailed in Appendix A.3. Heuristically, the class $\mathcal{G}_3(\beta, C)$ is close to the set of all $F_W$ that have at least $\beta$ uniformly bounded derivatives, and enjoy finite first absolute moment; and $\mathcal{G}_4(\beta, C)$ is identical, except for an analogous smoothness condition on $(f_W^{\mathrm{Ft}})'$ rather than on $f_W^{\mathrm{Ft}}$. The class $\mathcal{G}_5(C, k)$ is a set of $F_W$'s that have a bounded density and bounded moments of order $4(k+1)$, and $\mathcal{G}_6(\beta, C)$ is a set of $F_W$'s satisfying this moment assumption and for which the $j$th derivative of $f_W^{\mathrm{Ft}}(t)$ decreases at least as fast as $(1 + |t|)^{-\beta - j}$, where $0 \leq j \leq 2k+2$.

In the setting of quantile estimation, a small amount of smoothness in the vicinity of the true quantile seems necessary in order to perform the distribution inversion. The function classes $\mathcal{G}_7$ and $\mathcal{G}_8$ ensure this, together with, in the case of $\mathcal{G}_8$, constraining the true quantile not to be close to the origin. This is necessary in order to tease out the fact that, for quantile estimators as well as distribution estimators, convergence rates tend to be faster away from the origin than they are close to the origin.

3.2. *Upper bounds to convergence rates for $\widehat{F}_W$.* First we treat the case where $F_\delta$ is particularly "rough," in the sense that its characteristic function converges so slowly to zero in the tails that we may use the estimator $\widehat{F}_W$ with $h$ arbitrarily small; see Appendix A.1. Results (3.4) and (3.5), below, show that in this setting root-$n$ consistency is possible. A converse to (3.5) will be given in Theorem 3.4.

THEOREM 3.1. *Assume $\int |K| < \infty$ and $\int K = 1$. Then, for each $C_1, C_2 > 0$,*

$$(3.4) \quad \sup_{F_\delta \in \mathcal{F}_1(C_1)} \sup_{F_W \in \mathcal{G}_1(C_2)} \sup_{n \geq 1} \sup_{-\infty < x < \infty} nE\{\widehat{F}_W(x \mid 0) - F_W(x)\}^2 \leq 4C_1 C_2/\pi,$$

$$(3.5) \quad \sup_{F_\delta \in \mathcal{F}_2(C_1)} \sup_{F_W \in \mathcal{G}_2(C_2)} \sup_{n \geq 1} n \int_{-\infty}^{\infty} E\{\widehat{F}_W(x \mid 0) - F_W(x)\}^2 \, dx \leq 4(C_1 + C_2).$$

Next we treat cases where, in general, choosing a strictly positive value of $h$ can be advantageous. For definiteness, we choose $K$ so that $K^{\text{Ft}}$ is a compactly supported piece of a polynomial:

$$K^{\text{Ft}}(t) = (1 - t^r)^s 1(|t| \leq 1),$$

(3.6)

where $r \geq 2$ is an even integer, and $s \geq 1$ is an integer.

Such kernels are widely used in practice, where they have good numerical and theoretical performance; see Delaigle and Hall (2006). They satisfy the conditions imposed on $K$ in Theorem 3.1. More general kernels may be used, but they generally require stronger conditions defining the function classes.

Define $\ell_h = 1 + |\log h|$ if $\alpha = \frac{1}{2}$, and $\ell_h = 1$ otherwise. In Theorem 3.2 below, (3.7) and (3.8) give convergence rates uniformly in all $x$, and in $x$ not close to the origin, respectively. These rates are shown in Theorem 3.5 to be optimal in the respective cases, if $\alpha \neq \frac{1}{2}$. Result (3.9) gives the $L_2$ convergence rate.

THEOREM 3.2. *Assume $K$ satisfies (3.6) with $r > \beta + \frac{1}{2}$. Then, for each $C_1, C_2 > 0$, $0 \leq h \leq 1$ and $n \geq 1$,*

$$\sup_{F_\delta \in \mathcal{F}_3(\alpha,C_1)} \sup_{F_W \in \mathcal{G}_3(\beta,C_2)} \sup_{-\infty < x < \infty} E\{\widehat{F}_W(x \mid h) - F_W(x)\}^2$$

(3.7)
$$\leq B\{h^{2\beta} + n^{-1}(1 + h^{-(2\alpha-1)}\ell_h)\},$$

$$\sup_{F_\delta \in \mathcal{F}_3(\alpha,C_1)} \sup_{F_W \in \mathcal{G}_4(\beta,C_2)} \sup_{|x| > x_0} E\{\widehat{F}_W(x \mid h) - F_W(x)\}^2$$

(3.8)
$$\leq B\{x_0^{-2} h^{2\beta+2} + n^{-1}(1 + h^{-(2\alpha-1)}\ell_h)\},$$

$$\sup_{F_\delta \in \mathcal{F}_3(\alpha,C_1)} \sup_{F_W \in \mathcal{G}_3(\beta,C_2)} \int_{-\infty}^{\infty} E\{\widehat{F}_W(x \mid h) - F_W(x)\}^2 \, dx$$

(3.9)
$$\leq B\{h^{2\beta+1} + n^{-1}(1 + h^{-(2\alpha-1)}\ell_h)\},$$

*where, in each case, $B > 0$ depends only on $C_1$, $C_2$, $r$, $s$, $\alpha$ and $\beta$.*

Result (3.7), when $0 < \alpha < \frac{1}{2}$ and $\beta = 0$, is close to (3.4), although without an explicit formula for $B$ on the right-hand side. Note that when $0 < \alpha < \frac{1}{2}$ we may take $h = 0$ in (3.7).

To exhibit convergence rates, define $\ell = \log n$ if $\alpha = \frac{1}{2}$, and $\ell = 1$ if $\alpha > \frac{1}{2}$; put $h_1 = h_2 = h_3 = 0$ if $0 < \alpha < \frac{1}{2}$, and $h_j = C(\ell/n)^{1/(2\alpha+2\beta+j-2)}$ if $\alpha \geq \frac{1}{2}$, where $C > 0$; define $\rho_j = n^{-1}$ if $0 < \alpha < \frac{1}{2}$; and put $\rho_j = (\ell/n)^{(2\beta+j-1)/(2\alpha+2\beta+j-2)}$

if $\alpha \geq \frac{1}{2}$. The rates in (3.10), (3.11) and (3.12) below are obtained on taking $h = h_1$, $h_3$ and $h_2$ in (3.7), (3.8) and (3.9), respectively.

COROLLARY 3.3. *If $K$ satisfies (3.6) with $r > \beta + \frac{1}{2}$, and if $h_1, h_2, h_3$ are chosen as suggested above, then*

$$(3.10) \quad \sup_{F_\delta \in \mathcal{F}_3(\alpha, C_1)} \sup_{F_W \in \mathcal{G}_3(\beta, C_2)} \sup_{-\infty < x < \infty} E\{\widehat{F}_W(x \mid h_1) - F_W(x)\}^2 = O(\rho_1),$$

$$(3.11) \quad \sup_{F_\delta \in \mathcal{F}_3(\alpha, C_1)} \sup_{F_W \in \mathcal{G}_4(\beta, C_2)} \sup_{|x| > x_0} E\{\widehat{F}_W(x \mid h_3) - F_W(x)\}^2 = O(\rho_3),$$

$$(3.12) \quad \sup_{F_\delta \in \mathcal{F}_3(\alpha, C_1)} \sup_{F_W \in \mathcal{G}_3(\beta, C_2)} \int_{-\infty}^{\infty} E\{\widehat{F}_W(x \mid h_2) - F_W(x)\}^2 \, dx = O(\rho_2).$$

The rates $\rho_1$, $\rho_2$ and $\rho_3$ are in the order $\rho_3 < \rho_2 < \rho_1$. That is, mean-square convergence away from the origin is fastest, followed by convergence of mean integrated squared error, followed by mean-square convergence across the whole real line. The reason, as we shall show more explicitly in Theorem 3.5 and in Section 3.3 below, is that the estimator $\widehat{F}_W$ has difficulty in the neighborhood of the origin, and performs better outside that region. In approximate terms, its squared bias is of order $h^{2\beta}$ within radius $O(h)$ of the origin, and of order $h^{2(\beta+1)}$ a further distance away. Therefore, the squared-bias contribution to mean integrated squared error is of order $h(h^\beta)^2 = h^{2\beta+1}$.

Note, however, that this discussion is predicated on the assumption that the distribution of $\delta$ is symmetric, and the distribution of $W$ is in both $\mathcal{G}_3(\beta, C_2)$ and $\mathcal{G}_4(\beta, C_2)$. If, for example, $f_W^{\text{Ft}} = e^{itB_1 t}(1 + B_2|t|)^{-\beta}$, for real $B_1$ and $B_2 > 0$, then $F_W \in \mathcal{G}_3(\beta, C_2)$ for sufficiently large $C_2$, but $F_W$ does not lie in $\mathcal{G}_4(\beta, C_2)$ for any $C_2$ unless $B_1 = 0$. Therefore, a degree of centering at zero is being assumed. Of course, if we shift the center of the distribution of $W$ to $B_1$, then the results described in the previous paragraph continue to hold if we replace "the origin" by "$B_1$" throughout. This should be born in mind when interpreting discussion below.

These sizes of squared bias are reflected directly by the first terms on right-hand sides of (3.7)–(3.9). Moreover, as is suggested by the second terms there, and will be confirmed by the more detailed analysis in Section 3.5, error-about-the-mean properties of $\widehat{F}_W$ are very similar near the origin and away from the origin; their orders of magnitude do not alter.

3.3. *Lower bounds to convergence rates for $\widehat{F}_W$*. If $E|\delta| < \infty$, then (3.1) is equivalent to

$$(3.13) \quad \int_1^\infty t^{-2} f_\delta^{\text{Ft}}(t)^{-2} \, dt < \infty.$$

We know from (3.5) in Theorem 3.1 that, provided $E|\delta| < \infty$, (3.13) is sufficient for root-$n$ consistency of $\widehat{F}_W$, in the mean integrated squared error sense, uniformly over $F_W \in \mathcal{G}_2(C)$ for each fixed $C > 0$. Our next result shows that, under a mild additional assumption, (3.13) is also necessary for root-$n$ consistency.

THEOREM 3.4. *Let $\widehat{F}$ denote any measurable functional of the data. If $f_\delta$ is of bounded variation and $f_\delta^{\mathrm{Ft}}(t)$ is nonvanishing and eventually, for sufficiently large, positive $t$, monotone decreasing in $t$, and if, for some $C > 0$,*

$$(3.14) \qquad \sup_{F_W \in \mathcal{G}_2(C)} \int_{-\infty}^{\infty} E\{\widehat{F}(x) - F_W(x)\}^2 \, dx = O(n^{-1}),$$

*then (3.13) holds.*

Next we show that, despite the difficulty that $\widehat{F}_W$ can experience in a neighborhood of the origin, it converges there at the minimax-optimal rate. Likewise, it has optimal performance away from the origin. In particular, the convergence rates at (3.10) and (3.11) are both optimal. In view of what we have already learned, it is unsurprising that the rate of convergence of mean integrated squared error, in (3.12), is not optimal. Faster convergence rates can be achieved by using variable-bandwidth methods, where the bandwidth close to the origin is an order of magnitude smaller than that away from the origin.

Recall the definitions $\rho_1 = n^{-2\beta/(2\alpha+2\beta-1)}$ and $\rho_3 = n^{-(2\beta+2)/(2\alpha+2\beta+1)}$, appropriate for $\alpha > \frac{1}{2}$. Let $\mathcal{E}$ denote the class of measurable functionals of the data $X_1, \ldots, X_n$. Theorem 3.5, below, demonstrates optimality of the convergence rates given in (3.7) and (3.8) of Theorem 3.2.

THEOREM 3.5. *Let $F_\delta$ be a distribution for which $f_\delta^{\mathrm{Ft}}$ is real-valued and positive everywhere, and $|(f_\delta^{\mathrm{Ft}})^{(j)}(t)| \leq C_1(1+|t|)^{-\alpha-j}$ for all $t$ and for $j = 0, 1, 2$, where $C_1 > 0$. Then, provided $\alpha > \frac{1}{2}$, $x_1 \neq 0$, and $C_2 > 0$ is sufficiently large, there exists $C_3 > 0$ such that*

$$\inf_{\widehat{F} \in \mathcal{E}} \sup_{F_W \in \mathcal{G}_3(\beta, C_2)} E\{\widehat{F}(0) - F_W(0)\}^2 \geq C_3 \rho_1,$$

$$\inf_{\widehat{F} \in \mathcal{E}} \sup_{F_W \in \mathcal{G}_4(\beta, C_2)} E\{\widehat{F}(x_1) - F_W(x_1)\}^2 \geq C_3 \rho_3.$$

3.4. *Convergence rates of moment and quantile estimators.* Let $k \geq 0$ be an integer, and $q \in (2k, 2k+2)$. Define $\ell_{hq} = 1 + |\log h|$ if $\alpha = q + \frac{1}{2}$, and $\ell_{qh} = 1$ otherwise; and put $\rho_4 = n^{-(2\beta+2q)/(2\alpha+2\beta-1)}$. Result (3.16) below gives a convergence rate which, when $\alpha > q + \frac{1}{2}$ and $h = \mathrm{const.}\, n^{-1/(2\alpha+2\beta-1)}$, becomes identical to $O(\rho_4)$; and (3.17) shows that this rate is optimal. In that result we interpret $\bar{\nu}_q$ as a functional of $\widehat{F} \in \mathcal{E}$.

Theorem 3.6, below, is an analogue of Theorems 3.1 and 3.2 in the context of estimating the absolute moment $\nu_q$. It shows that root-$n$ consistency is possible, provided the distribution of $\delta$ is sufficiently rough; and it gives upper bounds to convergence rates in other cases.

THEOREM 3.6. *Let $k \geq 0$ be an integer, and let $2k < q < 2k+2$. Assume $K$ is given by (3.6), with, in the case of (3.15) below, $r > 2k+2$, and, for (3.16), $r > \max(\beta + q, 2k+2)$. Then, for each $C_1, C_2 > 0$ and for $0 \leq h \leq 1$,*

$$(3.15) \qquad \sup_{F_\delta \in \mathcal{F}_4(C_1, q)} \sup_{F_W \in \mathcal{G}_5(C_2, k)} E(\hat{\nu}_q - \nu_q)^2 \leq C n^{-1},$$

$$(3.16) \qquad \begin{aligned} \sup_{F_\delta \in \mathcal{F}_5(\alpha, C_1)} \sup_{F_W \in \mathcal{G}_6(\beta, C_2)} & E(\hat{\nu}_q - \nu_q)^2 \\ & \leq C\{h^{2\beta + 2q} + n^{-1}(1 + h^{-(2\alpha - 2q - 1)} \ell_{hq})\}, \end{aligned}$$

*where $C > 0$ depends only on $C_1, C_2$ and $q$. Furthermore, if $F_\delta \in \mathcal{F}_5(\alpha, C_1)$ and $\alpha > q + \frac{1}{2}$, then*

$$(3.17) \qquad \inf_{\widehat{F} \in \mathcal{E}} \sup_{F_W \in \mathcal{G}_6(\beta, C_2)} E(\bar{\nu}_q - \nu_q)^2 \geq C_3 \rho_4.$$

Theorem 3.4 has an analogue in this setting, asserting that if $E(\hat{\nu}_q - \nu_q)^2 = O(n^{-1})$ and $F_\delta$ satisfies mild additional assumptions, then the integrals at (A.7) (see Appendix A.2) converge.

To address the case of quantile estimation, recall from Section 2.3 the definitions of $\xi_u$ and $\hat{\xi}_u$, where $0 < u < 1$. Let $h_1$, $h_3$, $\rho_1$ and $\rho_3$ be as given immediately prior to Corollary 3.3, and let the function $g$ satisfy the conditions in the definition of $\mathcal{G}_u^*(C, g)$ in Section 3.1.

Results (3.18) and (3.19) below give upper bounds to rates of convergence for our quantile estimators when the quantile can lie anywhere, or is bounded away from the quantile for which $\xi_u = 0$, respectively. Results (3.20) and (3.21) give lower bounds, complementary to (3.18) and (3.19) respectively, in the case of general estimators.

THEOREM 3.7. *Assume that $\alpha > \frac{1}{2}$, with in addition $\alpha + 2\beta > 2$ and $\beta \geq 1$ in the case of (3.18) and (3.20), respectively; and suppose that $K$ is given by (3.6) with $r > \beta + \frac{1}{2}$. Then, for each $C_1, C_2 > 0$,*

$$(3.18) \quad \lim_{\lambda \to \infty} \limsup_{n \to \infty} \sup_{F_\delta \in \mathcal{F}_6(\alpha, C_1)} \sup_{F_W \in \mathcal{G}_7(\beta, C_2, u, g)} P\{|\hat{\xi}_u(h_1) - \xi_u| > \rho_1^{1/2} \lambda\} = 0,$$

$$(3.19) \quad \lim_{\lambda \to \infty} \limsup_{n \to \infty} \sup_{F_\delta \in \mathcal{F}_6(\alpha, C_1)} \sup_{F_W \in \mathcal{G}_8(\beta, C_2, u, g)} P\{|\hat{\xi}_u(h_3) - \xi_u| > \rho_3^{1/2} \lambda\} = 0,$$

$$(3.20) \quad \liminf_{\lambda \downarrow 0} \liminf_{n \to \infty} \inf_{\widehat{F} \in \mathcal{E}} \sup_{F_\delta \in \mathcal{F}_6(\alpha, C_1)} \sup_{F_W \in \mathcal{G}_7(\beta, C_2, u, g)} P(|\hat{\xi}_u - \xi_u| > \rho_1^{1/2} \lambda) > 0,$$

$$(3.21) \quad \liminf_{\lambda \downarrow 0} \liminf_{n \to \infty} \inf_{\widehat{F} \in \mathcal{E}} \sup_{F_\delta \in \mathcal{F}_6(\alpha, C_1)} \sup_{F_W \in \mathcal{G}_8(\beta, C_2, u, g)} P(|\hat{\xi}_u - \xi_u| > \rho_3^{1/2} \lambda) > 0.$$

3.5. *Limiting distributions.* Under conditions more restrictive than those imposed in Section 3.1, it is possible to obtain central limit theorems for $\widehat{F}_W$, exhibiting the convergence rates discussed in Section 3.2 and having explicitly-given biases and variances. The main features of these results are as follows: (a) The asymptotic variance equals a constant multiple of $n^{-1} h^{1-2\alpha}$, where the constant, $V(x)$, say, depends on $x$; (b) When $x = 0$, the asymptotic bias is a constant multiple, $B_1$, say, of $h^\beta$; and (c) When $x \neq 0$ the bias is asymptotic to $B_2(h,x) h^{\beta+1}$, where $B_2(h,x)$ is uniformly bounded as $h \downarrow 0$, and exceeds, in absolute value and for arbitrarily small $h$, a fixed constant as $h$ decreases. Formulae for $V(x)$, $B_1$ and $B_2(h,x)$ are given at (3.22) and (3.23).

To appreciate the relevance of these results, we interpret them in the context of (3.7) and (3.8). Excepting the case $\alpha = \frac{1}{2}$, there is a term of size $n^{-1} h^{1-2\alpha}$ on the right-hand sides of both those formulae. This term represents the main effect of variance, and is as indicated in (a) above. In (3.8), which is for the case of values $x$ that are bounded away from zero, there is a term of size $h^{2\beta+2}$ on the right-hand side. This represents the main effect of squared bias, and (c) above notes that its order of magnitude cannot be reduced. In (3.7), which includes the case $x = 0$, there is a term of size $h^{2\beta}$ on the right-hand side, and as (b) above observes, this too cannot be reduced.

Next we outline regularity conditions that give rise to these explicit expansions. Recall that most of the classes of distributions $F_\delta$ ask that $f_\delta^{\text{Ft}}$ decrease no faster than $t^{-\alpha}$ as $t$ increases. On the present occasion our main requirements are that $f_\delta^{\text{Ft}}$, and its first derivative, have an explicit expansion in inverses of polynomials up to a degree which strictly exceeds $2\alpha$, and that $f_W^{\text{Ft}}(t)$ behave to first order like a constant multiple of $t^{-\beta}$, with a remainder that is small enough to permit inversion of the characteristic function uniformly in $|x| > x_0$. These properties are captured by regularity conditions (A.8) and (A.9), given in Appendix A.4.

THEOREM 3.8. *Assume (A.8) and (A.9), that $f_X$ is bounded and continuous at $x$, and that the bandwidth satisfies $h = h(n) \to 0$ and $nh^{2\alpha-1} \to \infty$. Then, $\widehat{F}_W(x) - F_W(x)$ is asymptotically normally distributed with mean $B_1 h^\beta + o(h^\beta)$ or $B_2(h,x) h^{\beta+1} + o(h^{\beta+1})$, according as $x = 0$ or $x \neq 0$ respectively, where*

$$B_1 = -\frac{b}{2\pi} \sum_{j=1}^{s} (\,s\,) j \frac{(-1)^j}{rj - \beta} \qquad or$$

(3.22)
$$B_2(h,x) = -\frac{a\cos(x/h) + b\sin(x/h)}{2\pi x};$$

and with variance $n^{-1}h^{1-2\alpha}V(x)$, where

(3.23) $$V(x) = \pi^{-2}z^2 f_X(x) \int_0^\infty \left\{\int_0^1 \frac{\sin tu}{t}(1-t^r)^s t^\alpha\, dt\right\}^2 du.$$

In (3.22) and (3.23), the integers $r$ and $s$ are as at (3.6).

## 4. Numerical properties.

4.1. *Finite-sample performance of the distribution function, absolute moment and quantile estimators.* In this section we report the results of a simulation study illustrating the theoretical results and finite-sample behavior of the estimators of population features considered in Sections 2 and 3. We consider three distributions for $W$: (1) $W \sim N(0,1)$, (2) $W \sim \frac{1}{2}N(-3,1) + \frac{1}{2}N(2,1)$ and (3) $W \sim \text{Gamma}(2,1)$. Distributions (1), (2) and a variant of (3) were considered by Delaigle and Gijbels (2004a); see their #1, #3 and #2, respectively. Note that (1) gives a unimodal, symmetric (about 0) density, (2) a bimodal and two-sided density, and (3) a unimodal, one-sided density. Furthermore, the tails of the characteristic functions of (1) and (2) decay exponentially fast, while those of (3) decay at a polynomial rate.

For the error distribution, we consider the symmetrized Gamma $(\alpha,1)$ distributions [cf. (3.3)] with $\alpha = 2$ or $\alpha = 6$. For each combination of the target and error distributions, we consider two different sample sizes, $n = 100$ and $n = 800$, and a range of values of the smoothing parameter $h$, specifically $\{0.2, 0.4, \ldots, 2.0\}$. In the simulation study for this section we choose the kernel $K$ at (3.6), with $r = 4$ and $s = 2$. The number of simulation runs used in each case is 500.

The set of $x$ values is $\{-0.8, 0, 1.5\}$ for model (1), $\{-3.0, -0.5, 1.5\}$ for (2) and $x \in \{0.8, 1.5, 3.0\}$ for (3). Note that $x = 0$ is a common point of symmetry for $W$ and $\delta$ under model (1), while $x = -3.0$ and $x = -0.5$ are respectively one of the modes of $W$ and the mean (or median) of $W$ under model (2). The other $x$-values are chosen such that each set addresses both sides of the median of $W$. The optimal value of $h$ depends on the level, $u$, of the quantile, although the degree of sensitivity varies from one model to another. For all three models considered in our numerical work, the optimal $h$ lies in the interval $[0.4, 0.8]$ for the coarser error distribution with $\alpha = 2$, and in the interval $[1.0, 1.4]$ when $\alpha = 6$, reflecting the fact that a larger value of $h$ is more appropriate for error variables with a smoother distribution.

Figures 1 and 2 give the mean squared errors (MSEs) of the distribution function estimator $\widehat{F}_W(x|h)$ as a function of $h$ for different values of the argument $x$, under models (1) and (2) respectively. The graphs in the case of model (3) are close to those for model (1), provided $x = -0.8$, 0 and 1.5 are replaced by $x = 0.8$, 1.5 and 3.0, respectively.

The shape of the target distribution (i.e., unimodality versus bimodality) also seems to have an effect on the MSE curve, and hence, on the optimal value of the smoothing parameter, $h$. Interestingly, for $\alpha = 6$, the value $h = 1.0$ is the best choice, among those considered, for all the $x$'s and $n$'s

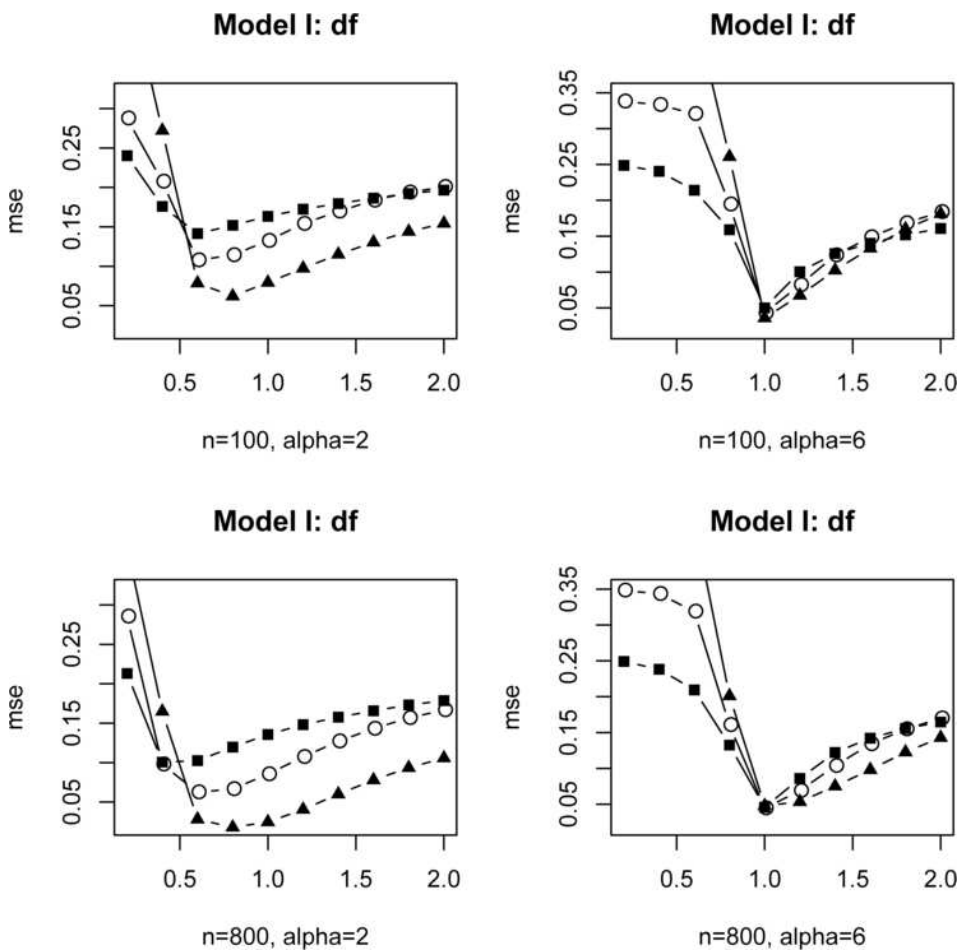

FIG. 1. *MSEs of the distribution function estimator $\hat{F}_W(x|h)$ under model (1), as a function of $h \in \{0.2, 0.4, \ldots, 2.0\}$, for $x \in \{-0.8, 0.0, 1.5\}$. In each panel, the MSE curves are marked with circles for $x = -0.8$, with squares for $x = 0.0$ and with triangles for $x = 1.5$. The error distribution is given by (3.3) with $\alpha \in \{2, 6\}$. The results are based on 500 simulation runs.*

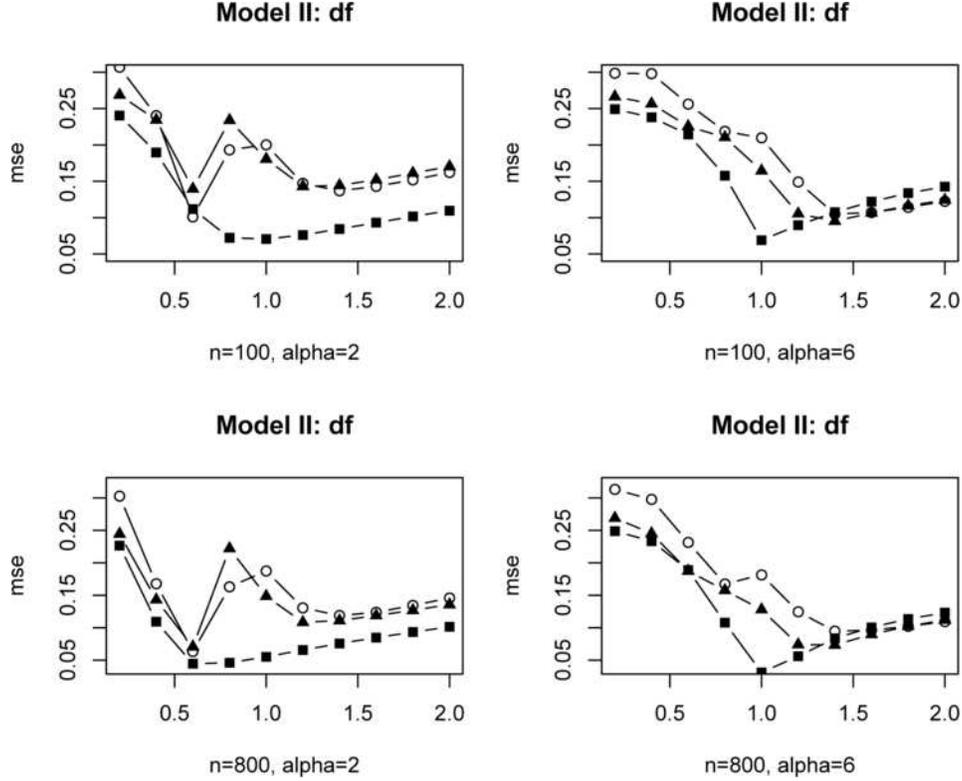

Fig. 2. *MSEs of the distribution function estimator $\hat{F}_W(x|h)$ under model (2), as a function of $h \in \{0.2, 0.4, \ldots, 2.0\}$, for $x \in \{-3.0, -0.5, 1.5\}$. In each panel the MSE curves are marked with circles for $x = -3.0$, with squares for $x = -0.5$ and with triangles for $x = 1.5$. The error distribution is given by (3.3) with $\alpha \in \{2, 6\}$. The results are based on 500 simulation runs.*

under models (1) and (3), and also for $x = -0.5$ and both the $n$'s under model (2).

Next, we consider the absolute moment estimator $\hat{\nu}_q(h)$ of Section 2.2, and the quantile estimator $\hat{\xi}_u(h)$ of Section 2.3. Figures 3 and 4 give the MSE functions of $\hat{\nu}_q(h)$ for $q \in \{0.5, 1, 1.5\}$, and of $\hat{\xi}_u(h)$ for $u \in \{0.4, 0.5, 0.7\}$ under model (1). To save space, we omit results for the other two models. The range of $h$ values is the same as before, except in the case of the absolute moment estimators where $h$ is restricted to a subset of $\{0.2, 0.4, \ldots, 1.4\}$. For $h \in \{1.6, 1.8, 2.0\}$, the MSE's of $\hat{\nu}_q(h)$ become too large (in 100s to 1000s, depending on the value of $q$), and hence, are omitted from the plot.

Note that the MSE functions for estimating the absolute moments are also nicely curved, in all cases attaining their minima, among the values of $h$ considered, at $h = 1.0$. However, the moment estimator $\hat{\nu}_q(h)$ seems to be very sensitive to oversmoothing, that is, to choice of too-high values of $h$. In comparison, the MSE functions of the quantile estimator $\hat{\xi}_u(h)$ are much more stable for under- and over-smoothing. Further, unlike the cases of moment and distribution function estimation, estimation of the two lower quantiles, $u = 0.4$ and $0.5$, seems to be less sensitive to smoothness of the error law; here the best performance is achieved when the values of $h$ are small. For the higher quantile, $u = 0.7$, the optimal $h$ shows dependence on the smoothness level of the error distribution, with larger $h$-values giving better performance.

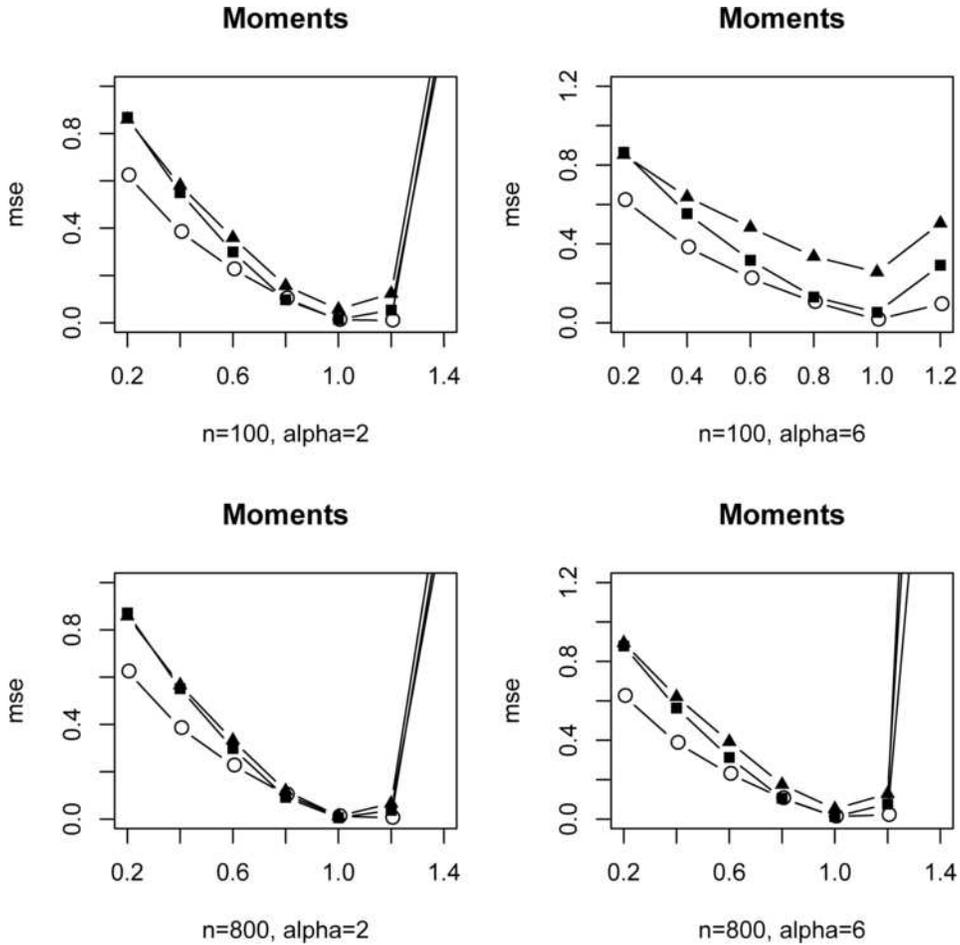

Fig. 3. *MSEs of the absolute moment estimator $\hat{\nu}_q(h)$ under model (1), as a function of $h$, for $q \in \{0.5, 1.0, 1.5\}$. In each panel the MSE curves are marked with circles for $q = 0.5$, with squares for $q = 1.0$ and with triangles for $q = 1.5$. The error distribution is given by (3.3) with $\alpha \in \{2, 6\}$. The results are based on 500 simulation runs.*

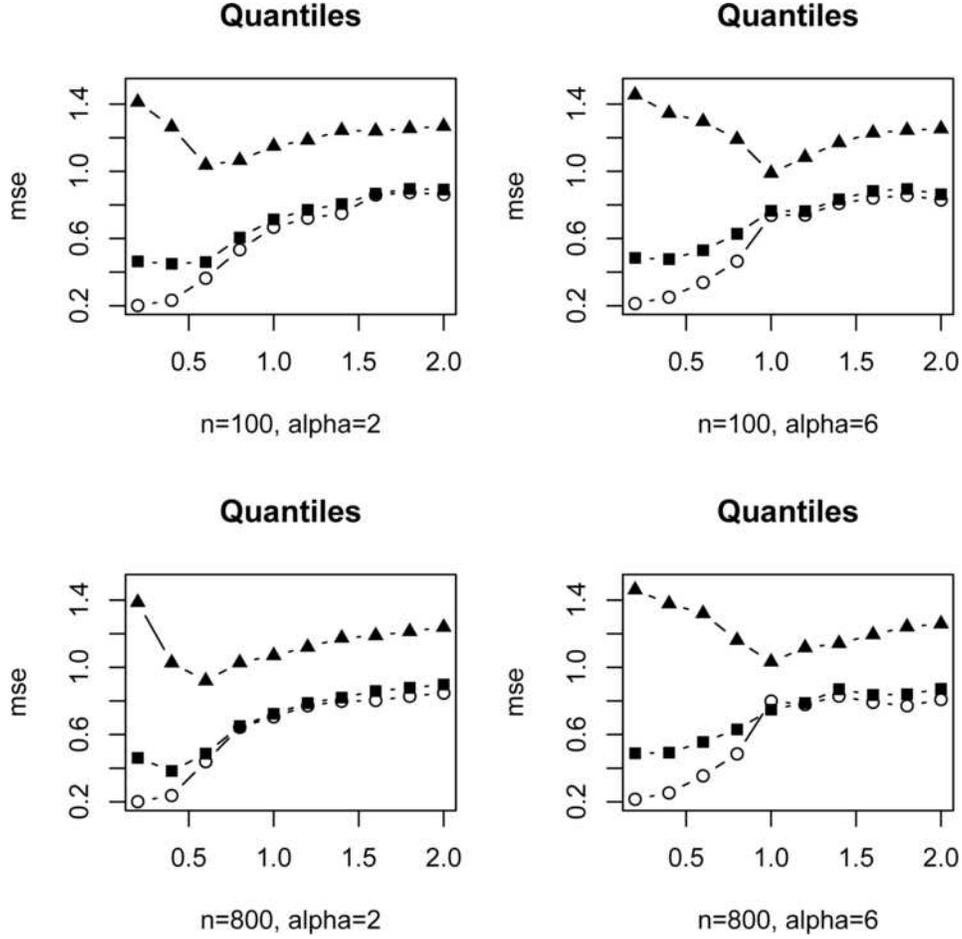

FIG. 4. *MSEs of the quantile estimator $\hat{\xi}_u(h)$ under model (1), as a function of $h \in \{0.2, 0.4, \ldots, 2.0\}$, for $u \in \{0.4, 0.5, 0.7\}$. In each panel the MSE curves are marked with circles for $u = 0.4$, with squares for $u = 0.5$ and with triangles for $u = 0.8$. The error distribution is given by (3.3) with $\alpha \in \{2, 6\}$. The results are based on 500 simulation runs.*

We next consider the effects of the argument on accuracy of distribution function and quantile estimation (cf. Theorems 3.2 and 3.7). Recall that under model (1), and under our choice of the symmetric error distribution, Theorem 3.2 asserts that the estimator $\widehat{F}_W(x)$ has a faster optimal rate of convergence at a nonzero $x$ compared to that at $x = 0$. Similar behavior is predicted by Theorem 3.7 for the quantile estimator $\hat{\xi}_u(h)$ at $u \neq 0.5$ and at $u = 0.5$. Figure 5 gives boxplots of the differences $\widehat{F}_W(x) - F_W(x)$ at $x = -0.8, 0, 1.5$ and $\hat{\xi}_u(h) - \xi_u$ at $u = 0.4, 0.5, 0.7$ under model (1) and $\alpha = 2$. The smoothing parameter for each value of the argument in $\widehat{F}_W(x)$ and $\hat{\xi}_u(h)$ is chosen to be the corresponding optimal value from Figures 1 and 4, respectively.

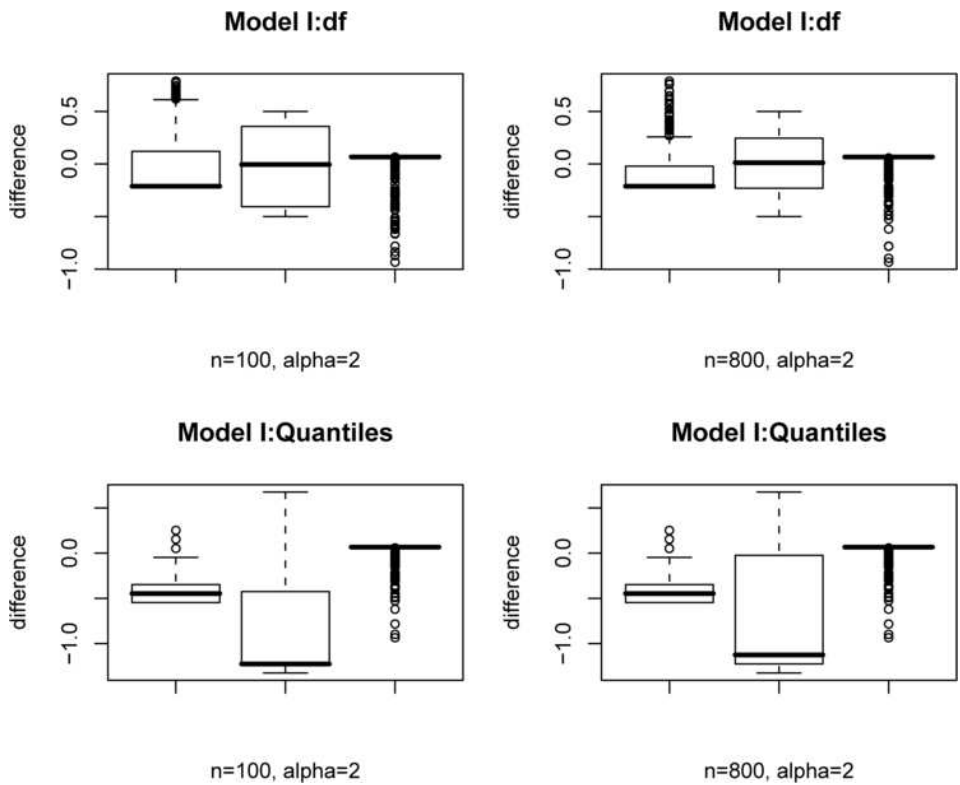

Fig. 5. *Box-plots of the deviations of the distribution function estimates and quantile estimates from their target values under model (1). In each of the top two panels the three box-plots correspond to the difference $\hat{F}_W(x|h) - F_W(x)$ for $x \in \{-0.8, 0.0, 1.5\}$, respectively, and in the lower two panels, to $\hat{\xi}_u(h) - \xi_u$ for $u \in \{0.4, 0.5, 0.8\}$, respectively. Here, the h-values are set at the respective optimal levels given in Figures 1 and 4. The results are based on 500 simulation runs.*

It is evident from Figure 5 that in the case of estimating the distribution function, the case $x = 0$ shows maximum variability around the target value at both sample sizes $n = 100$ and $n = 800$. The lower panels of Figure 5 show a similar pattern for the median estimator, $u = 0.5$, compared to the quantile estimators with $u = 0.4$ and $0.7$. The pronounced negative bias in the case of the median estimator reflects the difficulty of estimating the distribution function at the median [in the case of model (1)], and of estimating the median itself, relative to estimation at other places. See Theorem 3.7. However, the extent of the bias is greater than we had anticipated.

4.2. *Empirical choice of bandwidth.* We shall modify the "normal reference" approach, suggested by Delaigle and Gijbels (2004a) in the setting of density, rather than distribution, estimation. In particular, we shall temporarily take $f_W$ to be a normal $N(0, \sigma_W^2)$ density, with $\sigma_W^2 = \text{var}(W) = \text{var}(X) - \text{var}(\delta)$; and compute an estimator $\hat{\sigma}_W^2$ of $\sigma_W^2$ as the variance of the data $X_i$, minus the known variance of $\delta$. (The optimal bandwidth is invariant under changes to the location of $F_W$.)

To implement this approach, we shall use the following account of mean integrated squared error of the estimator $\widehat{F}_W(\cdot \mid h)$; see Appendix A.4 for regularity conditions.

THEOREM 4.1. *If (A.10) and (A.11) hold, then, as $n \to \infty$ and $h \to 0$,*

$$
\begin{aligned}
\int_{-\infty}^{\infty} & E\{\widehat{F}_W(x \mid h) - F_W(x)\}^2 \, dx \\
& = n^{-1} I(h) + B_W h^4 + o(n^{-1} h^{1-2\alpha} + h^4),
\end{aligned}
\tag{4.1}
$$

*where*

$$ 2\pi I(h) = \int t^{-2} \{1 - K^{\text{Ft}}(ht)/f_\delta^{\text{Ft}}(t)\}^2 \, dt, \qquad B_W = \tfrac{1}{4} \kappa_2^2 \int (f_W')^2 $$

*and $\kappa_2 = \int x^2 K(x) \, dx$, and $\alpha$ denotes the exponent of decay of $f_\delta^{\text{Ft}}(t)$.*

We may compute $I(h)$ by numerical integration. Alternatively, it can be approximated as $I(h) \sim A_\delta h^{1-2\alpha}$, where $A_\delta = C^2 \kappa / \pi$, $\kappa = \int_{t>0} t^{2\alpha - 2} K^{\text{Ft}}(t)^2 \, dt$, and the constants $C$ and $\pi$ are as in the asymptotic relation, $f_\delta^{\text{Ft}}(t) \sim C t^{-\alpha}$ as $t \to \infty$.

Bandwidth choice involves replacing $I(h)$ by its known value, for a particular $h$ (or using the approximation noted just above); replacing $B_W$ by its estimator, $\widehat{B}_W = \tfrac{1}{4} \kappa_2^2 (4\pi^{1/2} \hat{\sigma}_W^3)^{-1}$; and selecting $h$ by minimizing the resulting approximation to the sum of the first two terms on the right-hand side of (4.1). Note that, in the normal case, $R_W \equiv \int (f_W')^2 = (4\pi^{1/2} \sigma_W^3)^{-1}$ and so can be approximated by $\hat{R}_W \equiv (4\pi^{1/2} \hat{\sigma}_W^3)^{-1}$. In the results discussed below we used the exact value of $I(h)$.

We now report the results of a simulation study designed to investigate finite sample properties of this empirical bandwidth-selection procedure. We consider three distributions for $W$ as described above, namely, (1) $W \sim N(0,1)$, (2) $W \sim \tfrac{1}{2} N(-3,1) + \tfrac{1}{2} N(2,1)$, (3) $W \sim \text{Gamma}(2,1)$, and symmetrized Gamma $(\alpha, 1)$ distributions [cf. (3.3)] with $\alpha = 1$ and $\alpha = 5$ for the error distribution. Table 1 gives the theoretically optimal bandwidths obtained by minimizing the MISE in (4.1). Models (1) and (3) are seen to require almost identical amounts of smoothing, with model (2) needing a little more.

TABLE 1
*Theoretically optimal bandwidths for distribution function deconvolution*

| $(\alpha, n) =$ | (1,100) | (1,800) | (5,100) | (5,800) |
|---|---|---|---|---|
| Model (1) | 0.18 | 0.12 | 0.36 | 0.31 |
| Model (2) | 0.21 | 0.14 | 0.38 | 0.33 |
| Model (3) | 0.17 | 0.11 | 0.35 | 0.30 |

Following Delaigle and Gijbels (2004a), we used a one-step iteration method to compute $\hat{R}_W$; the optimal bandwidth estimator ($\hat{h}$, say) minimized the resulting estimated MISE function. Table 2 gives the bias and mean squared error of $\hat{h}$, based on 500 simulation runs. Numerical results not given here, for the sake of brevity, show that for all six combinations of the error and target distributions, the MSE of $\hat{h}$ decreased with sample size. Moreover, estimation is most accurate for distribution (1). This is likely due to use of the "normal reference" in the first step of the iteration. As expected, the performance of the method is better for the rougher error distribution ($\alpha = 1$), for all target distributions.

Next we consider the performance of the distribution-function estimators, using integrated squared error (ISE): $\int \{\hat{F}_W(x|\hat{h}) - F_W(x)\}^2 \, dx$. Table 3 gives values of the bias and the mean squared error of the ISE. Box-plots of scaled ISE values are given in Figure 6. In each case the scaling factor is the (theoretical) minimum of the MISE function. The distributions of the scaled ISE values behave as predicted by the theory for variations in sample size and in the smoothness of the error distribution. Further, from the box-plots it appears that, out of the three distributions considered here, the bimodal case (2) is the most difficult to recover.

## APPENDIX

**A.1. Definition of estimator $\widehat{F}(x \mid h)$.** The integral of $L$, the latter defined at (2.2), is given by $\int_{v \leq u} L(v) \, dv = L_1(hu)$, where, provided $K$ and $f_\delta$

TABLE 2
*The bias and the mean squared error (mse) of the estimated optimal bandwidths $\hat{h}$ based on 500 simulation runs. Here $(x)e(d)$ stands for $x \times 10^d$*

| $(\alpha, n) =$ | (1,100) | | (1,800) | | (5,100) | | (5,800) | |
|---|---|---|---|---|---|---|---|---|
| | bias | mse | bias | mse | bias | mse | bias | mse |
| Model (1) | 2.5e-2 | 1.1e-3 | 1.8e-2 | 3.8e-4 | 6.2e-2 | 2.8e-2 | 9.5e-3 | 6.6e-4 |
| Model (2) | 1.2e-1 | 1.5e-2 | 5.1e-2 | 2.7e-3 | 9.2e-2 | 9.0e-3 | 7.0e-2 | 5.1e-3 |
| Model (3) | 7.1e-2 | 5.4e-3 | 4.0e-2 | 1.6e-3 | 6.8e-2 | 1.3e-2 | 4.9e-2 | 2.5e-3 |

TABLE 3
*The bias and the mean squared error (mse) of the ISE of deconvolution distribution function estimators based on 500 simulation runs. Here $(x)e(d)$ stands for $x \times 10^d$*

| $(\alpha, n) =$ | (1,100) | | (1,800) | | (5,100) | | (5,800) | |
|---|---|---|---|---|---|---|---|---|
| | bias | mse | bias | mse | bias | mse | bias | mse |
| Model (1) | 4.4e-3 | 1.1e-4 | 5.3e-4 | 1.8e-06 | 1.4e-1 | 1.5e-1 | 2.0e-2 | 1.5e-3 |
| Model (2) | 1.6e-2 | 5.1e-4 | 1.8e-3 | 7.1e-06 | 3.4e-2 | 2.0e-3 | 7.1e-3 | 1.1e-4 |
| Model (3) | 5.9e-3 | 1.4e-4 | 7.5e-4 | 2.5e-06 | 5.1e-2 | 1.0e-2 | 3.5e-3 | 1.8e-4 |

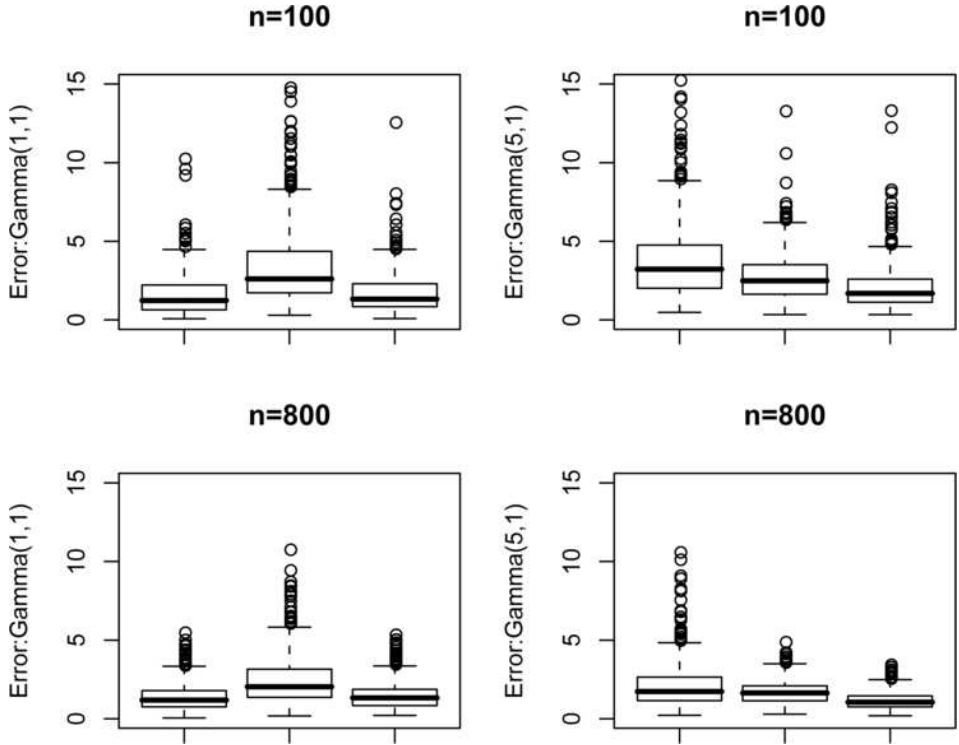

FIG. 6. *Box-plots of the scaled ISEs of the deconvolution distribution function estimators (scaled by the respective minimum MISEs). In each panel the three box-plots correspond to models (1), (2) and (3), respectively. The number of simulations was 500 in each case. Percentage of outliers falling outside the prescribed ranges of the boxplots are .05% for model (2) in the "top, left" panel, 5% for model (1) in the "top, right," none for the "bottom, left" and .05% for model (1) in the "bottom, right" panel, respectively.*

are both symmetric functions,

$$(A.1) \qquad L_1(u) = L_1(u \mid h) = \frac{1}{2} + \frac{1}{2\pi} \int_{-\infty}^{\infty} \frac{\sin t u}{t} \frac{K^{\text{Ft}}(ht)}{f_\delta^{\text{Ft}}(t)} \, dt.$$

Thus, by integrating $\hat{f}_W$, even if $f_W$ is not well defined, we obtain an estimator, $\widehat{F}_W$, of the distribution function, $F_W$, of $W$:

$$(A.2) \qquad \widehat{F}_W(x) = \widehat{F}_W(x \mid h) = \int_{-\infty}^{x} \hat{f}_W(u)\,du = \frac{1}{n}\sum_{j=1}^{n} L_1(x - X_j).$$

If $K^{\mathrm{Ft}}$ is compactly supported, and $f_\delta^{\mathrm{Ft}}$ does not vanish on the real line, then the integral at (A.1) is well defined and finite, provided $h \neq 0$. However, in view of Theorems 3.1 and 3.4, the case $h = 0$ is of particular interest. Since $\int K = 1$ then $K^{\mathrm{Ft}}(0) = 1$, and so it follows from (A.1) that

$$(A.3) \qquad L_1(u \mid 0) = \frac{1}{2} + \frac{1}{2\pi} \int_{-\infty}^{\infty} \frac{\sin tu}{t} \frac{1}{f_\delta^{\mathrm{Ft}}(t)}\,dt,$$

assuming that the integral on the right-hand side exists in the Riemann sense. An integration by parts argument shows that, for the integral in (A.3) to be Riemann convergent for each $u \neq 0$, it is sufficient that

(A.4)  $f_\delta^{\mathrm{Ft}}(t)$ is differentiable and $(d/dt)\{tf_\delta^{\mathrm{Ft}}(t)\}$ is integrable.

Reflecting (A.3), we take $L_1(u \mid 0) = \frac{1}{2}$ when $u = 0$.

The models for $f_\delta$ that are commonly used in practice are of Laplace type, and there

(A.5)  $|f_\delta^{\mathrm{Ft}}(t)|$ and $|t(d/dt)\{tf_\delta^{\mathrm{Ft}}(t)\}|$ are both bounded, both above and below, by constant multiples of $|t|^{-\alpha}$, as $|t|$ increases,

where $\alpha > 0$ is a parameter of the model. In this setting, (A.4) holds if and only if $\alpha < 1$, and then (A.3) also prevails. When (A.5) is true, the constraint $\alpha < 1$ is less constrictive than (3.1), which characterizes root-$n$ consistency of $\widehat{F}_W(\cdot \mid 0)$ for $F_W$; see Theorems 3.1 and 3.4. Indeed, if (A.5) holds, then (3.1) is true if and only if $\alpha < \frac{1}{2}$.

Therefore, for the sort of distribution of $\delta$ for which one might practically be interested in taking $h = 0$ in $\widehat{F}(x \mid h)$, one can expect the estimator

$$\widehat{F}_W(x \mid 0) = \frac{1}{n}\sum_{j=1}^{n} L_1(x - X_j \mid 0)$$

to be well defined and finite for each $x$. More generally, however, provided (3.1) obtains, the quantities

$$(A.6) \quad E\{\widehat{F}_W(x \mid 0) - F_W(x)\}^2 \quad \text{and} \quad \int_{-\infty}^{\infty} \{\widehat{F}_W(x \mid 0) - F_W(x)\}^2\,dx$$

are well defined and finite, either in their own right or as limits of their counterparts when $h > 0$, without considering models for which (A.4) holds. Existence in their own right follows from the fact that, assuming (3.1), (A.3) implicitly defines almost-everywhere a function $L_1(u \mid 0) - \frac{1}{2}$, which, by Parseval's theorem, is square-integrable. There are several ways of formally defining this function, for example, as an almost-everywhere limit of a subsequence of a sequence of Fourier inverses of compactly-supported approximations to the Fourier transform of $L_1(u \mid 0) - \frac{1}{2}$, or as an almost-everywhere limit along a subsequence, as $h \to \infty$, of $L_1(u \mid h) - \frac{1}{2}$.

Hence, it is appropriate to discuss the value of, and rate of convergence to zero of, both of the quantities at (A.6), without imposing conditions such as (A.5). Reflecting this point, in the formulation of Theorem 3.1 we do not require such assumptions.

**A.2. Classes of potential distributions of $\delta$.** Note particularly that all the function classes $\mathcal{F}_j$ include constraints which prevent $f_\delta^{\text{Ft}}$ from ever vanishing if the corresponding distribution lies in that class. Given $\alpha, C > 0$, write $\mathcal{F}_3(\alpha, C)$ for the class of continuous distributions $F_\delta$ for which $f_\delta^{\text{Ft}}$ is real-valued and positive, $\sup f_\delta \leq C$, $E|\delta| \leq C$ and $C f_\delta^{\text{Ft}}(t) \geq (1+|t|)^{-\alpha}$. (See the end of Appendix A.3 for interpretation of conditions on boundedness of densities.)

Given an integer $k \geq 0$, and $q \in (2k, 2k+2)$, let $\mathcal{F}_4(C, q)$ denote the class of $F_\delta$ for which $E(\delta^{4(k+1)}) \leq C$ and $f_\delta^{\text{Ft}}$ is real-valued and positive and satisfies

$$
\text{(A.7)} \quad \int_0^1 t^{-2(q+1)} \left\{ f_\delta^{\text{Ft}}(t) - \sum_{j=0}^k \frac{(-1)^j t^{2j}}{(2j)!} E(\delta^{2j}) \right\}^2 dt \\
+ \int_1^\infty t^{-2(q+1)} f_\delta^{\text{Ft}}(t)^{-2} dt \leq C.
$$

[The first part of (A.7) is essentially a moment condition.] Write $\mathcal{F}_5(\alpha, C)$ for the class of all $F_\delta \in \mathcal{F}_3(\alpha, C)$ for which $E(\delta^{4(k+1)}) \leq C$. Let $\mathcal{F}_6(\alpha, C)$ be the set of all $F_\delta \in \mathcal{F}_3(\alpha, C)$ for which $|(f_\delta^{\text{Ft}})^{(j)}(t)| \leq C(1+|t|)^{-\alpha-j}$ for $j = 0, 1, 2$. [Therefore, if $F_\delta \in \mathcal{F}_6(\alpha, C)$, then $|f_\delta^{\text{Ft}}|$ is bounded above and below by constant multiples of $(1+|t|)^{-\alpha}$.]

**A.3. Classes of potential distributions of $W$.** For $\beta \geq 0$ and $C > 0$, let $\mathcal{G}_3(\beta, C)$ be the class of $F_W$ for which $|f_W^{\text{Ft}}(t)| \leq C(1+|t|)^{-\beta}$ and $E|W| \leq C$, and let $\mathcal{G}_4(\beta, C)$ be the class of $F_W \in \mathcal{G}_3(\beta, C)$ such that

$$\sup_{u > 0}(1+u)^{\beta-k} \sup_{|x| > x_0} \left| \int_0^u e^{itx} t^k f_W^{\text{Ft}}(t)\, dt \right| \leq \frac{Ck}{(k_\beta - \beta)|x_0|}$$

for each $x_0 > 0$ and each integer $k > \beta$, where $k_\beta$ is the least such integer. (See two paragraphs below for interpretation of this constraint.) Given an integer $k \geq 0$, let $\mathcal{G}_5(C, k)$ be the class of $F_W$ for which $\sup f_W \leq C$ and $E(W^{4(k+1)}) \leq C$; and write $\mathcal{G}_6(\beta, C)$ for the class of $F_W$ such that $|(f_W^{\text{Ft}})^{(j)}(t)| \leq C(1+|t|)^{-\beta-j}$ for $0 \leq j \leq 2k+2$, and $E(W^{4(k+1)}) \leq C$.

Let $0 < u < 1$, let $g:[0,1] \to [0,\infty)$ be such that $g(x) \to 0$ as $x \downarrow 0$, and denote by $\mathcal{G}_u^*(C,g)$ the class of $F_W$ such that (a) $f_W$ exists and is strictly positive in $\mathcal{I}_u(C) = [\xi_u - C^{-1}, \xi_u + C^{-1}]$, (b) $f_W(\xi_u) + f_W(\xi_u)^{-1} \le C$, and (c) $|f_W(x) - f_W(y)| \le g(\eta)$ for all $\eta \in [0,1]$ and all $x,y \in \mathcal{I}_u(C)$ with $|x-y| \le \eta$. Write $\mathcal{G}_7(\beta,C,u,g)$ for the class of all $F_W \in \mathcal{G}_3(\beta,C) \cap \mathcal{G}_u^*(C,g)$, and $\mathcal{G}_8(\beta,C,u,g)$ for the class of all $F_W \in \mathcal{G}_4(\beta,C) \cap \mathcal{G}_u^*(C,g)$ for which $|\xi_u| > C^{-1}$.

Next we elucidate some of these function classes. If $C > 0$ is sufficiently large then $\mathcal{G}_3(\beta,C)$ contains $\phi_\beta$, defined at (3.3). To appreciate the sorts of distributions that are in $\mathcal{G}_4(\beta,C)$, note that in many instances where $F_W$ is centered at the origin and $F_W \in \mathcal{G}_3(\beta,C_1)$, it holds true that, for some $C_2 > 0$, $|(f_W^{\text{Ft}})'(t)| \le C_2(1+|t|)^{-\beta-1}$. Consider, for example, the case where $f_W^{\text{Ft}}(t) = (1+B|t|)^{-\beta}$ with $B > 0$. In such cases, an integration-by-parts argument shows that, for $u > 0$,

$$\left| x \int_0^u e^{itx} t^k f_W^{\text{Ft}}(t) \, dt \right| = \left| e^{iux} u^k f_W^{\text{Ft}}(u) - \int_0^u e^{itx} t^{k-1} \{ k f_W^{\text{Ft}}(t) + t (f_W^{\text{Ft}})'(t) \} \, dt \right|$$

$$\le \{ C_1 + (k_\beta - \beta)^{-1}(C_1 k + C_2) \} (1+u)^{k-\beta}.$$

Therefore, $F_W \in \mathcal{G}_4(\beta,C)$ if $C \ge 2C_1 + C_2$.

In the definitions of function classes $\mathcal{F}_j$ and $\mathcal{G}_j$, the conditions (a) $f_\delta \le C$ or (b) $f_W \le C$ are imposed only to ensure that (c) $\sup f_X \le C$. Property (c) holds if either (a) or (b) does, and so it is possible to switch the condition $\sup f_\delta \le C$, in the definition of a function class $\mathcal{F}_j$, to $\sup f_W \le C$, in the definition of $\mathcal{G}_k$, whenever both $f_\delta \in \mathcal{F}_j$ and $f_W \in \mathcal{G}_k$ are assumed. This feature allows variants of several of our theorems to be formulated easily; those variants will not be discussed explicitly.

**A.4. Regularity conditions for Theorems 3.8 and 4.1.** For Theorem 3.8 we take the kernel to be given by (3.6), with integers $r$ and $s$ as in that formula, and assume that

(A.8) the distribution of $\delta$ is symmetric about the origin, and for all $t \ge 0$,
$$f_\delta^{\text{Ft}}(t) = z^{-1}(1+t)^{-\alpha} + z_1(1+t)^{-\alpha_1} + \cdots + z_p(1+t)^{-\alpha_p} + A(t),$$
where $\frac{1}{2} < \alpha < \alpha_1 < \cdots < \alpha_{p+1}, \alpha_{p+1} > 2\alpha, z, z_1, \ldots, z_p$ are nonzero real numbers, and $|A^{(j)}(t)| \le \text{const.}(1+t)^{-\alpha_{p+1}-j}$ for $j = 0, 1$;

as $t \to \infty$, $f_W^{\text{Ft}}(t) = (a+ib)t^{-\beta} + o(t^{-\beta})$, where $a, b$ are real numbers, $\beta > 0$ and, for each $x_0 > 0$ and for $0 \le k \le rs - 1$,

(A.9)
$$(1+u)^{\beta-k} \sup_{|x|>x_0} \left| \int_0^u e^{itx} t^k \{ f_W^{\text{Ft}}(t) - (a+ib)(1+t)^{-\beta} \} \, dt \right| \to 0$$

as $u \to \infty$.

In Theorem 4.1 we assume that

(A.10) $K$ is symmetric and satisfies $\int K = 1$ and $\kappa_2 \equiv \int x^2 K(x)\,dx \neq 0$, and $K^{\mathrm{Ft}}$ is compactly supported;

(A.11) $E|\delta| < \infty$; for all $t, f_\delta^{\mathrm{Ft}}(t) \neq 0$; for some $\alpha > \frac{1}{2}$ and $C > 0$, $f_\delta^{\mathrm{Ft}}(t) \sim Ct^{-\alpha}$ as $t \to \infty$; and for some $\beta > \frac{3}{2}$ and $C_1 > 0$, $F_W \in \mathcal{G}_3(\beta, C_1)$.

The condition $\beta > \frac{3}{2}$ in (A.11) implies that $\int (f_W')^2 < \infty$.

**Acknowledgments.** We are grateful to two referees, an associate editor and an editor, whose constructive comments encouraged this revision in which technical details are pushed a distance into the background.

## REFERENCES

BOOTH, J. G. and HALL, P. (1993). Bootstrap confidence regions for functional relationships in errors-in-variables models. *Ann. Statist.* **21** 1780–1791. MR1245768

BUTUCEA, C. (2004). Deconvolution of supersmooth densities with smooth noise. *Canad. J. Statist.* **32** 181–192. MR2064400

BUTUCEA, C. and TSYBAKOV, A. B. (2008). Sharp optimality for density deconvolution with dominating bias. *Theory Probab. Appl.* To appear.

CARROLL, R.J. and HALL, P. (1988). Optimal rates of convergence for deconvolving a density. *J. Amer. Statist. Assoc.* **83** 1184–1186. MR0997599

CORDY, C. and THOMAS, D. R. (1997). Deconvolution of a distribution function. *J. Amer. Statist. Assoc.* **92** 1459–1465. MR1615256

CUI, H. (2005). Asymptotics of mean transformation estimators with errors in variables model. *J. Syst. Sci. Complex.* **18** 446–455. MR2172109

DELAIGLE, A. and GIJBELS, I. (2002). Estimation of integrated squared density derivatives from a contaminated sample. *J. Roy. Statist. Soc. Ser. B* **64** 869–886. MR1979392

DELAIGLE, A. and GIJBELS, I. (2004a). Practical bandwidth selection in deconvolution kernel density estimation. *Comput. Statist. Data Anal.* **45** 249–267. MR2045631

DELAIGLE, A. and GIJBELS, I. (2004b). Bootstrap bandwidth selection in kernel density estimation from a contaminated sample. *Ann. Inst. Statist. Math.* **56** 19–47. MR2053727

DELAIGLE, A. and HALL, P. (2006). On optimal kernel choice for deconvolution. *Statist. Probab. Lett.* **76** 1594–1602. MR2248846

DEVROYE, L. (1989). Consistent deconvolution in density estimation. *Canad. J. Statist.* **17** 235-239. MR1033106

DIGGLE, P. J. and HALL, P. (1993). A Fourier approach to nonparametric deconvolution of a density estimate. *J. Roy. Statist. Soc. Ser. B* **55** 523–531. MR1224414

FAN, J. (1991a). On the optimal rates of convergence for nonparametric deconvolution problems. *Ann. Statist.* **19** 1257–1272. MR1126324

FAN, J. (1991b). Global behavior of deconvolution kernel estimates. *Statist. Sinica* **1** 541–551. MR1130132

FAN, J. (1993). Adaptively local one-dimensional subproblems with application to a deconvolution problem. *Ann. Statist.* **21** 600–610. MR1232507

Fan, J. and Koo, J.-Y. (2002). Wavelet deconvolution. *IEEE Trans. Inform. Theory* **48** 734–747. MR1889978

Groeneboom, P. and Jongbloed, G. (2003). Density estimation in the uniform deconvolution model. *Statist. Neerlandica* **57** 136–157. MR2035863

Groeneboom, P. and Wellner, J. (1992). *Information Bounds and Nonparametric Maximum Likelihood Estimation*. Birkhäuser, Basel. MR1180321

Hesse, C. H. (1995). Distribution function estimation from noisy observations. *Publ. Inst. Stat. Paris Sud* **39** 21–35.

Hesse, C. H. (1999). Data-driven deconvolution. *J. Nonparametr. Statist.* **10** 343–373. MR1717098

Hesse, C. H. and Meister, A. (2004). Optimal iterative density deconvolution. *J. Nonparametr. Statist.* **16** 879–900. MR2094745

Ioannides, D. A. and Papanastassiou, D. P. (2001). Estimating the distribution function of a stationary process involving measurement errors. *Statist. Inference Stoch. Process.* **4** 181–198. MR1856173

Jongbloed, G. (1998). Exponential deconvolution: Two asymptotically equivalent estimators. *Statist. Neerlandica* **52** 6–17. MR1615570

Koo, J.-A. (1999). Logspline deconvolution in Besov space. *Scand. J. Statist.* **26** 73–86. MR1685303

Neumann, M. H. (1997). On the effect of estimating the error density in nonparametric deconvolution. *J. Nonparametr. Statist.* **7** 307–330. MR1460203

Pensky, M. (2002). Density deconvolution based on wavelets with bounded supports. *Statist. Probab. Lett.* **56** 261–269. MR1892987

Pensky, M. and Vidakovic, B. (1999). Adaptive wavelet estimator for nonparametric density deconvolution. *Ann. Statist.* **27** 2033–2053. MR1765627

Qin, H.-Z. and Feng, S.-Y. (2003). Deconvolution kernel estimator for mean transformation with ordinary smooth error. *Statist. Probab. Lett.* **61** 337–346. MR1959070

Stefanski, L. A. and Carroll, R. J. (1990). Deconvoluting kernel density estimators. *Statistics* **21** 169–184. MR1054861

Stone, C. J. (1982). Optimal global rates of convergence for nonparametric regression. *Ann. Statist.* **10** 1040–1053. MR0673642

van de Geer, S. (1995). Asymptotic normality in mixture models. *ESAIM Probab. Statist.* **1** 17–33. MR1382516

van Es, B., Spreij, P. and van Zanten, H. (2003). Nonparametric volatility density estimation. *Bernoulli* **9** 451–465. MR1997492

Zhang, C. H. (1990). Fourier methods for estimating mixing densities and distributions. *Ann. Statist.* **18** 806–830. MR1056338

Department of Mathematics and Statistics  
University of Melbourne  
Melbourne, VIC 3010  
Australia

Department of Statistics  
Texas A&M University  
College Station, Texas 77843-3143  
USA  
E-mail: snlahiri@tamu.edu